\documentclass[11pt]{article}

\usepackage[top=1in, bottom=1in, left=0.85 in, right=0.85in]{geometry}

\usepackage{graphicx,color,float,epstopdf}
\usepackage[linesnumbered,ruled]{algorithm2e}
\usepackage{multirow,array}
\usepackage{pst-blur,pstricks-add}
\usepackage{latexsym, bm, bbm, amsthm, amsmath,amssymb}
\usepackage{cases}
\usepackage{epsfig}
\usepackage{epstopdf}
\usepackage{contour,soul,color}
\usepackage[colorlinks=true]{hyperref}
\hypersetup{urlcolor=blue, citecolor=blue,linkcolor=blue}

\usepackage{amsfonts}
\usepackage{mathtools}
\usepackage{cite}
\usepackage[numbers]{natbib}
\setcitestyle{open={},close={}}

\def\epi{{\rm epi}\,}

\def\ri{{\rm ri}\,}
\def\sgn_A{{\rm sgn}_A\,}
\def\dom{{\rm dom}\,}

\def\gph{{\rm gph}\,}
\def\dom{{\rm dom}\,}
\def\dist{{\rm dist}\,}

\def\rank{{\rm rank}\,}
\def\KKT{{\rm KKT}\,}

\newcommand{\dd}{d\kern -0.16em l}

\newtheorem{theorem}{Theorem}[section]

\newtheorem{Lemma}{Lemma}[section]
\newtheorem{corollary}{Corollary}[section]
\newtheorem{definition}{Definition}[section]

\newtheorem{example}{Example}[section]

\begin{document}

\vskip -1.5cm
\title{Critical Multipliers in Semidefinite Programming
\footnote{
Supported by the National Natural Science Foundation of China under project grant No.11571059, No.11731013 and No. 91330206.
}}
\author{Tianyu Zhang\footnote{School of Mathematical Sciences, Dalian University of Technology,  China. (tyzhang@mail.dlut.edu.cn).}
\,\,and\,
 Liwei Zhang\footnote{School of Mathematical Sciences, Dalian University of Technology,  China. (lwzhang@dlut.edu.cn).} 
}
\vspace{2mm}
\date{\today}
\maketitle

\begin{abstract}  
It was proved in [\cite{izmailov2014newton}] that the existence of a noncritical multiplier for a (smooth) nonlinear programming problem is equivalent to an error bound condition for the Karush-Kuhn-Thcker (KKT) system without any assumptions. This paper investigates whether this result still holds true for a (smooth) nonlinear semidefinite programming (SDP) problem. We first  introduce the notion of critical and noncritical multipliers for a SDP problem and obtain their  complete characterizations  in terms of the problem data. We prove for the SDP problem, the noncriticality property can be derived from the error bound condition for  the KKT system without any assumptions, and this fact is revealed by some simple examples. Besides we give an appropriate second-order sufficient optimality condition characterizing noncriticality explicitly.  We propose a set of  assumptions from which  the error bound condition for the KKT system can be derived from the noncriticality property. Finally we establish a new error bound for $x$-part, which is expressed by both perturbation and the multiplier estimation.
\end{abstract}

\section{Introduction}
Consider an  optimization problem of the following form
\begin{equation}\label{1.1}\left\{
\begin{array}{ll}
    \min\limits_{x\in X} & f(x) \\
  {\rm s.t.} &  G(x)\in K,
 \end{array}\right.
\end{equation}
where $f:X\rightarrow \mathbb{R}$, $G : X\rightarrow Y$, $X$ and $Y$ are two finite dimensional real vector spaces, and $K$ is a closed convex set in $Y$.
We say that a local error bound holds in a neighborhood of a  particular  KKT point $(\bar{x},\bar{y})$  if there exist a function $\sigma: X\times Y\rightarrow[0,\infty)$ and a constant $c>0$ such that
\begin{equation}\label{1.2}
 \sigma(x,y)\geq c \,\dist\big((x,y),\bar{x}\times M(\bar{x})\big),
\end{equation}  for all $(x,y)$ in a neighborhood of $(\bar{x},\bar{y})$, where $M(\bar{x})$ is the set
of all the multipliers at the point $\bar{x}$.

When $f$ and $G$ are continuously differentiable  and
$K$ is a polyhedral set, Problem (\ref{1.1}) is  a conventional nonlinear programing (NLP). It is well known that for NLP the local error bound plays a critical role in the convergence analysis for numerical algorithms when the  KKT system has nonisolated solutions. For example, if the error bound condition holds, algorithms can be constructed  which converge locally superlinearly to a KKT point [\cite{izmailov2014newton}, Chapter 7]. Therefore the characterization of the error bound condition is an important issue for the study of NLP problem.  Izmailov and solodov defined critical
multipliers and noncritical conterparts, they showed the existence of the local error bound (\ref{1.2}) around the KKT point $(\bar{x},\bar{y})$ is equivalent to requiring that $\bar{y}$ is noncritical [\cite{izmailov2014newton}, Chapter 1].

Recently Mordukhovich and Sarabi [\cite{mordukhovich2017critical}] consider the following problem
\begin{equation}\label{1.3}
\min f(x)=\varphi_0(x)+\theta(\Phi(x)), x\in\mathbb{R}^n,
\end{equation}
where  $\Phi:\mathbb{R}^n\rightarrow\mathbb{R}^m$ and $\Psi:\mathbb{R}^n\times\mathbb{R}^m\rightarrow\mathbb{R}^l$ are sufficiently smooth while the convex function $\theta:\mathbb{R}^m\rightarrow\bar{\mathbb{R}}:=(-\infty,\infty]$ is  extended-real-valued but piecewise linear. This is just the case of Problem (\ref{1.1}) where $f(x):=\varphi_0(x)+\theta(\Phi(x))$, $G(x):=\Phi(x)$, and $K:={\rm \dom}\theta:=\{z\in\mathbb{R}^m|\theta(z)<\infty\}$.  They proposed  critical and noncritical multipliers for the KKT system of (\ref{1.3}), and established the equivalent relation between noncritical multipliers and the local error bound (\ref{1.2}).

When $K$ is a nonpolyhedral set, critical and noncritical multipliers are not easy to be characterized. Sun, Cui and Toh [\cite{Cui2016On}] provided sufficient conditions to guarantee the metric subregularity of KKT mappings for solving
linearly constrained convex semidefinite programming with multiple solutions. When the set $K$ in (\ref{1.1}) belongs to the class of $C^2$- cone reducible sets, Ding, Sun and Zhang [\cite{ding2017characterization}] showed that under the Robinson constaint qualification, the KKT solution mapping is robustly isolated calm if and only if both the strict robinson constraint qualification and second order sufficient condition hold. Cui, Ding and Zhao [\cite{cui2017quadratic}] provided two types of sufficient conditions for ensuring the quadratic
growth conditions of a class of constrained convex symmetric and non-symmetric matrix optimization problems regularized by nonsmooth spectral functions. These sufficient conditions are derived via the study of the $C^2$-cone reducibility of spectral functions and the metric subregularity of their subdifferentials, respectively.

Mordukhovich and Sarabi [\cite{mordukhovich2017critical}] introduced the notions of critical and noncritical multipliers for KKT systems of a class of composite optimization problems which is a extension of the notions by Izmailove and Solodov [\cite{Izmailov2013A}] developed for
classical KKT systems for NLPs.

   When $K=S^p_{+}$, the cone of positively semidefinite matrices in $S^p$, how can we define critical and noncritical multipliers, is the result for  the equivalence of  the existence of noncritical multipliers and the local error bound still valid?  These questions are quite interesting
and challenging.  The purpose of this paper is to try answering  these questions.

The remaining parts of this paper are organized as follows. In the next section, we introduce some definitions and preliminary results required in the subsequent analysis from variational analysis. In section 3, we define critical and noncritical multipliers for KKT systems of SDPs, establish equivalent description of them and specify them for noncritical multipliers for particular KKT systems. We also characterize noncritical multipliers for KKT systems for SDPs via some error bounds, present a new second-order sufficient condition
(SOSC$\bar{y}$) which ensures the noncriticality of multiplier $\bar{y}$. Section 4 is mainly devoted to characterizing another kind of error bound. We also proof (SOSC$\bar{y}$) is a sufficient condition for this kind of error bound. We conclude this paper in section 5.
\section{Preliminaries}
In this section we recall some basic notions and definitions in variational analysis which are extensively used in this paper.
For a set valued mapping $S:\mathcal {X}\rightrightarrows\mathcal {Y}$ with $\mathcal {X}$, $\mathcal {Y}$ being finite dimensional Hilbert spaces, the upper limite is the set
\begin{equation*}
\limsup\limits_{x\rightarrow\bar{x}}S(x):=\{y\in\mathcal {Y}:\exists x^{k}\rightarrow \bar{x}, ~\exists y^{k}\rightarrow \bar{y}, ~y^k\in S(x^k)\},
\end{equation*}
while the inner limit is the set
\begin{equation*}
\liminf\limits_{x\rightarrow\bar{x}}S(x):=\{y\in\mathcal {Y}:\forall x^{k}\rightarrow \bar{x},~\exists N_1\in \mathcal {N}_\infty , ~\exists y^{k}\overset{N_1}{\rightarrow} \bar{y}, ~y^k\in S(x^k),~k\in{N_1}\},
\end{equation*}
with $\mathcal {N}_\infty:=\{N\subseteq \mathcal {N} :\mathcal {N}\backslash N ~is ~finite\}$($N$ is the natural numbers).
Based on the definition of inner and outer limits of set-valued mappings, normal cone and tangent cone can be calculated. $\mathcal {Z}$ is a finite dimensional Hilbert space, and $\Omega\in \mathcal {Z}$ is locally closed around $\bar{x}\in\Omega$. The regular normal cone to $\Omega$~at $\bar{x}$ is defined by
\begin{equation}\label{1}
 \hat{N}(\bar{x};\Omega):=\left\{v\in\mathcal {Z}^\ast:\limsup\limits_{x\overset{\Omega} {\rightarrow}\bar{x}}\frac{\langle v,x-\bar{x}\rangle}{\|x-\bar{x}\|}\leq 0\right\},
\end{equation}
and normal cone by
\begin{equation}\label{2}
 {N}(\bar{x};\Omega):=\limsup\limits_{x\overset{\Omega} {\rightarrow}\bar{x}}\hat{N}(\bar{x};\Omega),
\end{equation}
where the symbol $x\overset{\Omega} {\rightarrow}\bar{x}$
means that $x {\rightarrow}\bar{x}$ with $x\in\Omega$. If
$\Omega$ is a convex set, then (\ref{1}) and (\ref{2}) are equal to the classical normal cone of convex analysis.
Define the normal cone and inner tangent cone to $\Omega$ at $\bar{x}$ by
\begin{equation}\label{3}
  T(\bar{x};\Omega):=\limsup\limits_{t\searrow 0}\frac{\Omega-\bar{x}}{t}
\end{equation}and
\begin{equation}\label{4}
  T^{i}(\bar{x};\Omega):=\liminf\limits_{t\searrow 0}\frac{\Omega-\bar{x}}{t}
\end{equation}
respectively.
If $h:\mathcal {Y}\rightarrow \bar{\mathbb{R}}$ is an extend-valued function on a finite Hilbert space, define the basic subdifferential of $h$ at $\bar{z}\in \dom~h$ by
\begin{equation}\label{5}
 \partial h(\bar{z}):=\{y\in\mathcal {Y}:(y,-1)\in N((\bar{z},h(\bar{z});\epi h))\}.\end{equation}
The indicator function for set $\Omega$ is $\delta_ \Omega=\delta(z;\Omega):=0$ for $z\in\Omega$ and $\delta(z;\Omega):=\infty$ otherwise, and
we have
\begin{equation}\label{6}
N(\bar{z};\Omega)=\partial \delta(\bar{z};\Omega),~ \bar{z}\in\Omega.
\end{equation}
The regular coderivative, the limiting coderivative and the graphical derivative to $S$ at $(\bar{x},\bar{y})\in \gph S $ are defined by
\begin{equation}\label{7}
\hat{D}^ \ast S(\bar{x},\bar{y})(v):=\{u\in\mathcal {X}:(u,-v)\in \hat{N}((\bar{x},\bar{y});\gph S)\},~v\in\mathcal {Y},
\end{equation}
\begin{equation}\label{8}
{D}^ \ast S(\bar{x},\bar{y})(v):=\{u\in\mathcal {X}:(u,-v)\in {N}((\bar{x},\bar{y});\gph S)\},~v\in\mathcal {Y},
\end{equation}
\begin{equation}\label{9}
 {D} S(\bar{x},\bar{y})(u):=\{v\in\mathcal {Y}:(u,v)\in {T}((\bar{x},\bar{y});\gph S)\},~u\in\mathcal {X},
\end{equation}
respectively.
\begin{Lemma}\label{Lemma1}[\cite{wu2014mathematical}, Theorem 3.1]
For any $(X,Y)\in \gph N(\cdot ; S^p_{+})$, the inner tangent cone and tangent cone to $\gph N(\cdot ;S^p_{+})$
coincide and they both have the expression
\begin{equation}\label{10}
\begin{array}{lll}
T^{i}((X,Y);\gph N(\cdot; S^p_{+}))&=&T((X,Y);\gph N(\cdot; S^p_{+}))\\&=& \{(H_1, H_2)\in{S^p\times S^p}:\Pi ^{\prime}_{S^p_{+}}(X+Y;H_1+H_2)=H_1\}.
\end{array}
\end{equation}
\end{Lemma}
For any matrix $A\in S^p$, let $A$ have the following eigenvalue decomposition
\begin{equation}\label{11}
A=P\left[
\begin{matrix}
  \Lambda_\alpha &  &  \\
   & 0_\beta &  \\
   & & \Lambda_\gamma
\end{matrix}
\right]P^T,
\end{equation}
where $\alpha:=\{i:\lambda_i(A)>0\}$, $\beta:=\{i:\lambda_i(A)=0\}$, $\gamma:=\{i:\lambda_i(A)<0\}$. Let $(X,Y)\in \gph N(\cdot,S^p_+)$ and $A=X+Y$ have the eigenvalue decomposition (\ref{11}) then the following results can be get easily
\begin{equation}\label{12}
X=P\left[
\begin{matrix}
 \Lambda_\alpha &  &  \\
   & 0_\beta &  \\
   &  & 0_\gamma
\end{matrix}
\right]P^T,~Y=P\left[
\begin{matrix}
 0_\alpha &  &  \\
   & 0_\beta &  \\
   &  &  \Lambda_\gamma
\end{matrix}
\right]P^T.
\end{equation}
Define the matrix $\sum\in S^P$ with entries
\begin{equation}\label{13}
 \Sigma_{ij}:=\frac{\max\{\lambda_{i}(A),0\}-\max\{\lambda_{j}(A),0\}}
  {\lambda_{i}(A)-\lambda_{j}(A)},i,j=1,\ldots,n,
\end{equation}
where $0/0$ is defined to be $1$. Denote ${P}=[{P_\alpha} {P_\beta} {P_\gamma}]$ with ${P_\alpha}\in\mathbb{R}^{p\times\alpha}$, ${P_\beta}\in\mathbb{R}^{p\times\beta}$,
${P_\gamma}\in\mathbb{R}^{p\times\gamma}$. Then from [\cite{Cui2016On}] we have
\begin{equation*}
\left\{
\begin{array}{lll}
  T(X;S^p_+)& = & \{H\in S^p:[{P_\beta} P_\gamma]^TH[{P_\beta} P_\gamma]\succeq 0\}, \\
  T(Y;S^p_-)& = &\{H\in S^p:[{P_\beta} P_\gamma]^TH[{P_\alpha} P_\beta]\preceq 0\}, \\
  N(X;S^p_+)& = & \{H\in S^p:[{P_\beta} P_\gamma]^TH[{P_\beta} P_\gamma]\preceq 0,~ {P_\alpha}^TH{P}=0\}, \\
  N(Y;S^p_-)& = & \{H\in S^p:[{P_\alpha} P_\beta]^TH[{P_\alpha} P_\beta]\succeq 0,~ {P_\gamma}^TH{P}=0\}.
\end{array}\right.
\end{equation*}
For the convenience of the latter discussions, we  denote the critical cone of $S^p_+$ at $X$ associated with $Y$ as
\begin{equation*}
C_{S^p_+}(X,Y):=T(X;S^p_+)\cap Y^\perp=\{H\in S^p:{P_\gamma}^TH[P_\beta P_\gamma]=0,~ P_\beta^THP_\beta\succeq0\},
\end{equation*}
and the critical cone of $S^p_-$ at $Y$ associated with $X$ as
\begin{equation*}
C_{S^p_-}(Y,X):=T(Y;S^p_-)\cap X^\perp=\{H\in S^p:{P_\alpha}^TH[P_\alpha P_\beta]=0,~ P_\beta^THP_\beta\preceq0\}.
\end{equation*}
\begin{Lemma}\label{Lemma2}[\cite{wu2014mathematical}, Corallary3.1]
For any $(X,Y)\in \gph N(\cdot;S^p_+)$, let $A=X+Y$ have the eigenvalue decomposition (\ref{11}). Then
\begin{equation}\label{14}
T((X,Y);\gph N(\cdot; S^p_{+}))=\left\{(H_1, H_2)\in{S^p\times S^p}: 
\begin{array}{ll}
& (\tilde{H_1})_{\beta\alpha}=0,(\tilde{H_1})_{\gamma\gamma}=0,
\\&(\tilde{H_2})_{\alpha\alpha}=0,(\tilde{H_2})_{\alpha\beta}=0, \\
   &(\Sigma_{\alpha\gamma}-E_{\alpha\gamma})\circ(\tilde{H_1})_
   {\alpha\gamma}+\Sigma_{\alpha\gamma}\circ(\tilde{H_2})_
   {\alpha\gamma}=0, \\
   & S_{+}^{|\beta|}\ni(\tilde{H_1})_{\beta\beta}\perp
   (\tilde{H_2})_{\beta\beta}\in S_{-}^{|\beta|}
\end{array}\right\},
\end{equation}
where $\tilde{H_1}=P^T H_1 P$, $\tilde{H_2}=P^T H_2 P$.
\end{Lemma}
\begin{Lemma}(Hoffman's Lemma)\label{Lemma3}
Let $X$ and $Y$ be Banach spaces, and let $A :X\rightarrow Y$ be a linear continuous mapping with closed range (i.e., its range $AX$ is a closed subspace of $Y$).  Given $x_i^\ast\in X^\ast$, $i=1,\cdots,p$, consider the multifunction
\begin{equation}\label{39}
\Upsilon(y,b):=\{x\in X: Ax=y, \langle x_i^\ast,x\rangle\leq b_i, i=1,\cdots,p\}.
\end{equation}
Then there exist a constant $k>0$, depending on $A$ and $x^\star_i, i=1,\cdots,p$, such that for any $x\in X$ and any $(y,b)\in \dom(\Upsilon)$,
\begin{equation}\label{40}
\dist(x,\Upsilon(y,b))\leq k(\|Ax-y\|+\sum\limits^{p}_{i=1}[\langle x_i^\ast,x\rangle-b_i]_+).
\end{equation}
\end{Lemma}
In our subsequent discussion, like [\cite{Cui2016On}] we need the concept of bounded linear regularity of a collection of closed convex sets, see, e.g., [\cite{Bauschke_1996}, Definition 5.6].
\begin{definition}\label{definition2}
Let $D_1, D_2, \cdots, D_m\subseteq\mathcal{X}$ be closed convex sets for some positive integer $m$. Suppose that $D:=D_1\cap D_2\cap\cdots,\cap ~D_m$ is non-empty. The collection $\{D_1, D_2, \cdots, D_m\}$ is said to be boundedly linearly regular if for every bounded set $\mathcal{B}\subseteq\mathcal{X}$, there exists a constant $\kappa >0$ such that
\begin{equation}\label{87}
\dist(x, D)\leq\kappa \max\{\dist(x, D_1), \cdots, \dist(x, D_m)\},~\forall x\in\mathbb{B}.
\end{equation}
\end{definition}
A sufficient condition to guarantee the property of bounded linear regularity was established in [\cite{bauschke1999strong}, Corollary 3].
\begin{Lemma}\label{Lemma5}
Let $D_1, D_2, \cdots, D_m\subseteq\mathcal{X}$ be closed convex sets for some positive integer $m$. Suppose that $D_1,  D_2, \cdots, D_r$ are polyhedrals for some $r\in\{0,1,\cdots, m\}$. Then a sufficient condition for $D_1, D_2, \cdots, D_m$ to be boundedly linearly regular is
\begin{equation}\label{86}
\bigcap\limits_{i=1,2,\cdots,r}
D_i\cap\bigcap\limits_{i=r+1,\cdots,m}\ri(D_i)\neq\emptyset.
\end{equation}
\end{Lemma}
\begin{Lemma}\label{Lemma7}[\cite{Cui2016On}, Proposition 3.2]
Let $X\in S^p_+$ and $Y\in N(X;S^p_+)$. Suppose that $X$ and $Y$ have the eigenvalue decompositions as in (\ref{12}). Then it holds that :
\begin{description}
\setlength{\itemsep}{0pt}
\setlength{\parsep}{0pt}
\setlength{\parskip}{0pt}
\item[(a)]$N(X;S^p_+)$ is a polyhedral set if and only if $|\alpha|\geq p-1$;\\
\item[(b)]$0\in Y+\ri(N(X;S^p_+))$ if and only if $|\beta|=0$, i.e., $\rank(X)+\rank(Y)=p$.
\end{description}
\end{Lemma}
Similar to [\cite{Cui2016On}, Proposition 3.3, Proposition 3.4], we can easily obtain the following conclusions.
\begin{Lemma}\label{Lemma8}
Let $X\in S^p_+$ and $Y\in N(X;S^p_+)$. Then $N(\cdot;S^p_+)$ is metrically subregular at $X$ for $Y$.
\end{Lemma}

\begin{Lemma}\label{Lemma6}
Let $\bar{X}\in S^p_+$ and $\bar{Y}\in N(X;S^p_+)$. Suppose that
$\bar{X}$ and $\bar{Y}$ have the same eigenvalue decompositions as $X$ and $Y$ respectively in
(\ref{12}). Then for all $(X,Y)\in S^n\times S^n$ with $Y\in N(X;S^p_+)$  sufficiently close 
to $(\bar{x},\bar{y})\in S^n\times S^n$, we have
\begin{equation}\label{88}
\left\{
\begin{array}{lll}
  \tilde{X}_{\alpha\alpha}=\Lambda_\alpha+O(\|\Delta X\|),& \tilde{X}_{\alpha\beta}=O(\|\Delta X\|), &\tilde{X}_{\alpha\gamma}=O(min\{\|\Delta X\|, \|\Delta Y\|\})\\
   \tilde{X}_{\beta\beta}=O(\|\Delta X\|), &\tilde{X}_{\beta\gamma}=O(\|\Delta X\|\|\Delta Y\|),  & \tilde{X}_{\gamma\gamma}=O(\|\Delta X\|\|\Delta Y\|), \\
  \tilde{Y}_{\alpha\alpha}=O(\|\Delta X\|\|\Delta Y\|), &  \tilde{Y}_{\alpha\beta}=O(\|\Delta X\|\|\Delta Y\|), &
  \tilde{Y}_{\alpha\gamma}=O(min\{\|\Delta X\|, \|\Delta Y\|\})\\
 \tilde{Y}_{\beta\beta}=O(\|\Delta Y\|), &\tilde{Y}_{\beta\gamma}=O(\|\Delta Y\|),&
 \tilde{Y}_{\gamma\gamma}=\Lambda_\gamma+O(\|\Delta Y\|),\\
\end{array}\right.
\end{equation}
\begin{equation}\label{89}
\tilde{Y}_{\alpha\gamma}+\Lambda_\alpha^{-1}
\tilde{X}_{\alpha\gamma}\Lambda_\gamma=O(\|\Delta X\|\|\Delta Y\|),
\end{equation}
where $\Delta X=X-\bar{X}$, $\Delta Y=Y-\bar{Y}$.
\end{Lemma}
\section{Noncritical multipliers for semidefinite programming }
Consider the following SDP problem
\begin{equation}\label{15}
\begin{array}{ll}
  \min & f(x) \\
  {\rm s.t.} &  G(x)\in S^{p}_{+},\\
\end{array}
\end{equation}
where $f : \mathbb{R}^n\rightarrow \mathbb{R}$ and $G : \mathcal {X}\rightarrow S^p$ are twice continuously differentiable functions.
The Lagrangian function of (\ref{15}) is defined by
\begin{equation}\label{16}
L(x,y):=f(x)+\langle y,G(x)\rangle, ~(x,y)\in\mathbb{R}^n \times S^p.
\end{equation}
Let $\bar{x}\in\mathbb{R}^n$ be a feasible point  to (\ref{15}). We define the set of Lagrange multipliers associated with $\bar{x}$ by
\begin{equation}\label{17}
\Lambda(\bar{x}):=\{y\in S^p:\nabla_xL(\bar{x},y)=0,~y\in N(G(\bar{x});S^p_+)\}.
\end{equation}
Suppose that $\bar{x}$ is a feasible solution to problem
(\ref{15}). The critical cone $C(\bar{x})$ of (\ref{15})
at $\bar{x}$ is defined by
\begin{equation}\label{77}
C(\bar{x}):=\{d\in\mathbb{R}^n:G^{\prime}(\bar{x})d\in T(G(\bar{x});{S^p_+}),~ f^{\prime}(\bar{x})d\leq0\}.
\end{equation}
If $\bar{x}$ is a stationary point of problem (\ref{15}) and $\bar{y}\in \Lambda(\bar{x})$, then
\begin{equation}\label{78}
\begin{array}{lll}
C(\bar{x})&=&\{d\in\mathbb{R}^n:G^{\prime}(\bar{x})d\in T(G(\bar{x});{S^p_+}),~ f^{\prime}(\bar{x})d=0\}\\
&=&\{d\in\mathbb{R}^n:G^{\prime}(\bar{x})d\in C_{S^p_+}(G(\bar{x}),\bar{y})\}.
\end{array}
\end{equation}
We say that the Robinson constraint qualification (RCQ) for problem (\ref{15}) holds at a feasible point $\bar{x}$ if
\begin{equation}\label{80}
G^{\prime}(\bar{x})\mathbb{R}^n+T(G(\bar{x});S^p_+)=S^p.
\end{equation}
It is well known that the RCQ holds at a locally optimal solution $\bar{x}\in\mathbb{R}^n$ if and only if $\Lambda(\bar{x})$
is a nonempty, convex, and compact subset of $S^p$.\\
The SRCQ is said to hold for problem (\ref{15}) at $\bar{x}$ with respect to $\bar{y}\in\Lambda(\bar{x})\neq\emptyset$ if
\begin{equation}\label{81}
G^{\prime}(\bar{x})\mathbb{R}^n+T(G(\bar{x});S^p_+)\cap \bar{y}^\perp=S^p.
\end{equation}
The set of Lagrange multipliers $\Lambda(\bar{x})$ is a singleton if the SRCQ holds.
Let $\bar{x}$ be a stationary of (\ref{15}), and $\bar{y}$ be an associated Lagrange multiplier. Assume SRCQ is satisfied at $\bar{x}$. Local optimality of $\bar{x}$ implies the second order necessary condition (SONC) in the following form
\begin{equation}\label{76}
\langle d,\nabla^2_{xx}L(\bar{x},\bar{y})d\rangle-\sigma(\bar{y},
T^2_{S^p_+}(G(\bar{x}),G^{\prime}(\bar{x})d))\geq0,~\forall d\in C(\bar{x}).
\end{equation}
Recall that the  second order sufficient condition (SOSC) which ensures that $\bar{x}$ is a locally optimal solution to (\ref{15}) is defined by
\begin{equation}\label{75}
\sup\limits_{\bar{y}\in \Lambda(\bar{x})}\{\langle d,\nabla^2_{xx}L(\bar{x},\bar{y})d\rangle-\sigma(\bar{y},
T^2_{S^p_+}(G(\bar{x}),G^{\prime}(\bar{x})d))\}>0,~\forall d\in C(\bar{x})\backslash\{0\}.
\end{equation}
\begin{definition}(SOSC$\bar{y}$)
We say  that the second order sufficient condition for problem (\ref{15}) satisfied at $\bar{x}$  respect to $\bar{y}\in \Lambda(\bar{x})$ is
\begin{equation}\label{yy}
\langle d,\nabla^2_{xx}L(\bar{x},\bar{y})d\rangle-\sigma(\bar{y},
T^2_{S^p_+}(G(\bar{x}),G^{\prime}(\bar{x})d))>0,~\forall d\in C(\bar{x})\backslash\{0\}.
\end{equation}
\end{definition}
\begin{definition}(critical and noncritical multipliers)\label{definition1}
A Lagrange multiplier $\bar{y}\in\Lambda(\bar{x})$ is called critical multiplier for (\ref{15}) if there exists $0\neq\xi\in\mathbb{R}^n$ satisfying the following generalized KKT system
\begin{equation}\label{18}
0\in \nabla^2_{xx}L(\bar{x},\bar{y})\xi+G^\prime(\bar{x})^\ast(D\partial\delta(\cdot;S^p_+))
(G(\bar{x}),\bar{y})(G^\prime(\bar{x})\xi).
\end{equation}
The multiplier $\bar{y}\in \Lambda(\bar{x})$ is noncritical for (\ref{15}) otherwise, i.e., when the generalized equation (\ref{18}) admits only the trivial solution $\xi=0$. \end{definition}
For NLP it has been proved that superlinear convergence results for dual stabilization methods
under assumptions that rely on the dual estimates being close enough to a noncritical Lagrange multiplier, see [\cite{robinson2015comments}].

In the following, we consider the characterization of noncritical multipliers.

\begin{theorem}\label{theorem1}Let $(\bar{x}, \bar{y})$ is a \KKT point of (\ref{15}).
Then $\bar{y}$ is noncritical for (\ref{15}) if and only if
we can derive $\xi=0$ from \begin{equation}\label{98}
\left\{\begin{array}{l}
\nabla^2_{xx}L(\bar{x},\bar{y})\xi+G^\prime(\bar{x})^\ast\eta=0, \\
 G^\prime(\bar{x})\xi-\Pi^{\prime}_{S^p_{+}}
 (G(\bar{x})+\bar{y};G^\prime(\bar{x})\xi+\eta)=0,
 ~(G^{\prime}(\bar{x})\xi,\eta)\in\mathbb{R}^n\times S^p.
\end{array}\right.
\end{equation}
\end{theorem}
\noindent{\bf Proof.}
Relation (\ref{18}) can be reformulated as
\begin{equation}\label{99}
\left\{\begin{array}{l}
\nabla^2_{xx}L(\bar{x},\bar{y})\xi+G^\prime(\bar{x})^\ast\eta=0, \\\eta\in(D\partial\delta(\cdot;S^p_+))
(G(\bar{x}),\bar{y})(G^\prime(\bar{x})\xi).
\end{array}\right.
\end{equation}
From (\ref{6}) we have
\begin{equation}\label{31}
\partial\delta(\cdot;S^p_+)=N(\cdot;S^p_+).
\end{equation}
In view of the definition of the graphical derivative (see (\ref{9})), we have
\begin{equation}\label{20}
\eta\in(D\partial\delta(\cdot;S^p_+))
(G(\bar{x}),\bar{y})(G^\prime(\bar{x})\xi)
\end{equation}
which is equivalent to
\begin{equation}\label{21}
(G^\prime(\bar{x})\xi,\eta)\in T((G(\bar{x}),\bar{y});\gph N(\cdot;S^p_+)).
\end{equation}
Note that $\bar{y}\in N(G(\bar{x});S^p_+)\Leftrightarrow
 S^p_+\ni G(\bar{x})\perp\bar{y}\in S^p_-$, so $G(\bar{x})$ and $\bar{y}$ have the same eigenvalue decomposition. Applying Lemma \ref{Lemma1}, we obtain
\begin{equation}\label{22}
\Pi ^{\prime}_{S^p_{+}}(G(\bar{x})+\bar{y};G^\prime(\bar{x})\xi+\eta)=
G^\prime(\bar{x})\xi.
\end{equation}This proof is completed.
\hfill $\Box$
\begin{corollary}\label{corollary1}
Assume $(\bar{x},\bar{y})$ is a \KKT point of (\ref{15}).
$\bar{y}$ is critical for (\ref{15}) if and only if the following system
\begin{equation}\label{100}
\left\{\begin{array}{l}
\nabla^2_{xx}L(\bar{x},\bar{y})\xi+G^\prime(\bar{x})^\ast\eta=0, \\
 G^\prime(\bar{x})\xi-\Pi^{\prime}_{S^p_{+}}
 (G(\bar{x})+\bar{y};G^\prime(\bar{x})\xi+\eta)=0,
 ~(G^{\prime}(\bar{x})\xi,\eta)\in\mathbb{R}^n\times S^p
\end{array}\right.
\end{equation}
admits a solution pair $(\xi,\eta)\in\mathbb{R}^n\times S^p$ with $\xi\neq 0.$
\end{corollary}

The following result provides us a useful criterion for charactering critical multipliers.

\begin{corollary}(equivalent description of critical multipliers).\label{corollary2}
Let $(\bar{x},\bar{y})$ be a \KKT point of (\ref{15}). Then $\bar{y}$ is critical for (\ref{15}) if and only if the following system
\begin{eqnarray}\label{19}
\nabla^2_{xx}L(\bar{x},\bar{y})\xi+G^\prime(\bar{x})^\ast\eta=0,\\
(P^T G^\prime(\bar{x})\xi P)_{\beta\gamma}=0, (P^T G^\prime(\bar{x})\xi P)_{\gamma\gamma}=0,\\
(P^T \eta P)_{\alpha\alpha}=0, (P^T \eta P)_{\alpha\beta}=0,\\
(\Sigma_{\alpha\gamma}-E_{\alpha\gamma})\circ(P^T G^\prime(\bar{x})\xi P)_
   {\alpha\gamma}+\Sigma_{\alpha\gamma}\circ(P^T \eta P)_
   {\alpha\gamma}=0,\\
S_{+}^{|\beta|}\ni(P^T G^\prime(\bar{x})\xi P)_{\beta\beta}\perp
   (P^T \eta P)_{\beta\beta}\in S_{-}^{|\beta|}\},\\
   (G^\prime(\bar{x})\xi,\eta)\in S^p\times S^p
\end{eqnarray}
admits a solution pair $(\xi,\eta)\in\mathbb{R}^n\times S^p$ with $\xi\neq 0.$
\end{corollary}
\noindent{\bf Proof.}
From the proof of the Theorem \ref{theorem1}, we obtain that (\ref{18}) can be reformulated as
\begin{equation}
\left\{\begin{array}{l}
\nabla^2_{xx}L(\bar{x},\bar{y})\xi+G^\prime(\bar{x})^\ast\eta=0, \\(G^\prime(\bar{x})\xi,\eta)\in T((G(\bar{x}),\bar{y});\gph N(\cdot;S^p_+)).
\end{array}\right.
\end{equation}
Then applying Lemma \ref{Lemma2},
the results can be derived.\hfill $\Box$

We next show noncritical multipliers of SDP is compatible with the NLP case.
\begin{example}\label{example1}
Consider \rm{SDP} (\ref{15}) with $G(x)={\rm{Diag}}(g_1(x),\cdots,g_p(x))$, where $g_i: \mathbb{R}^n\rightarrow\mathbb{R},~ i\in\{1,\cdots,p\}$. We denote $g(x):=(g_1(x),\cdots,g_p(x))^T$, $\bar{Y}:={\rm{Diag}}(\bar{\lambda}_1,\cdots,\bar{\lambda}_p)\in\Lambda(\bar{x})$, and $\bar{y}:=(\bar{\lambda}_1,\cdots,\bar{\lambda}_p)^T$. Suppose without loss of generality that $\bar{x}$ satisfies $g_1(\bar{x})=\cdots=g_p(\bar{x})=0$. Denote the index subsets by
\begin{equation}\label{67}
I_-:=\{i\in\{1,\cdots,p\}: \bar{\lambda}_i<0\}~~and~~
I_0:=\{i\in\{1,\cdots,p\}: \bar{\lambda}_i=0\}.
\end{equation}
From Theorem \ref{theorem1} we have that $\bar{Y}$ is a noncritical multiplier of $\bar{x}$ if and only if
\begin{equation}\label{68}
\left\{
\begin{array}{l}
\nabla^2_{xx}L(\bar{x},\bar{y})\xi+g^\prime(\bar{x})^\ast\eta=0,\\
\eta=(\eta_1,\cdots,\eta_p)\in\mathbb{R}^p_-,\\ \eta_i\nabla g_i(\bar{x})^T\xi=0~if~i\in I_0,\\
\nabla g_i(\bar{x})^T\xi=0~if~i\in I_-,
\nabla g_i(\bar{x})^T\xi\geq0~if~i\in I_0
\end{array} \right.\Longrightarrow\xi=0.
\end{equation}
The above relations also illustrate that $\bar{y}$ is a noncritical multiplier at $\bar{x}$ for the following  inequality constrained  NLP:
\begin{equation}\label{69}
\begin{array}{ll}
  \min & f(x) \\
  {\rm s.t.} &  g(x)\geq0.\\
\end{array}
\end{equation}
\end{example}
Now, we turn to explore the relation between noncriticality property and error bound condition. The KKT system of (\ref{15}) is defined by:
\begin{equation}\label{24}
0\in\left[
\begin{array}{c}
 \nabla_xL(x,y) \\
  -G(x)
\end{array}
\right]+\left[
\begin{array}{c}
  0\\
  N(y;S^p_{-})
\end{array}
\right].
\end{equation}
The canonical perturbation of system (\ref{24}) is defined by
\begin{equation}\label{25}
\left[\begin{array}{c}
  p_1 \\ p_2
\end{array}
\right]
\in\left[
\begin{array}{c}
 \nabla_xL(x,y) \\
  -G(x)
\end{array}
\right]+\left[
\begin{array}{c}
  0\\
  N(y;S^p_{-})
\end{array}
\right]
\end{equation}
with the canonical parameter pair $(p_1,p_2)\in \mathbb{R}^n\times S^p$. Consider the set-valued mapping
 $\Phi:\mathbb{R}^n\times S^p\rightrightarrows\mathbb{R}^n\times S^p$ associated with (\ref{25}) by
\begin{equation}\label{26}
\Phi(x,y):=
\left[
\begin{array}{c}
 \nabla_xL(x,y) \\
  -G(x)
\end{array}
\right]+\left[
\begin{array}{c}
  0\\
  N(y;S^p_{-})
\end{array}
\right]
\end{equation}
and then define the solution map $S_{KKT}:\mathbb{R}^n\times S^p\rightrightarrows\mathbb{R}^n\times S^p$ to (\ref{25})
as the inverse to (\ref{26}) by
\begin{equation}\label{27}
S_{\rm KKT}(p_1,p_2):=\{(x,y)\in\mathbb{R}^n\times S^p:(p_1,p_2)\in\Phi(x,y)\}.
\end{equation}
For each given $(p_1,p_2)\in\mathbb{R}^n\times S^p$, it is easily verified that the  the solution set of the KKT system
(\ref{25}) can be rewritten as
\begin{equation}\label{28}
S_{\rm KKT}(p_1,p_2)=\{(x,z-\Pi_{S^p_{+}}(z))\in\mathbb{R}^n\times S^p:\Psi(x,z)=(p_1,-p_2)\},
\end{equation}
where $\Psi : \mathbb{R}^n\times S^p\rightarrow\mathbb{R}^n\times S^p$ is Robinson's normal mapping defined by
\begin{equation}\label{29}
\Psi(x,z)=\left[
\begin{array}{c}
         \Psi_1(x,z) \\
         \Psi_2(x,z)
       \end{array}\right]=
\left[\begin{array}{c}
          \nabla f(x)+G^{\prime}(x)^\ast (z-\Pi_{S^p_+}(z)) \\
            G(x)-\Pi_{S^p_+}(z)
          \end{array}\right],~(x,z)\in\mathbb{R}^n\times S^p.
\end{equation}
Let $(\bar{x},\bar{y})$ be a solution to the KKT system (\ref{25}) with $(p_1,p_2)=(0,0)$. Denote $\bar{z}:=G(\bar{x})+\bar{y}$. Then $\Pi_{S^p_+}$ is globally Lipschitz continuous (with modulus 1) and $G$ is locally Lipschitz continuous.
\begin{theorem}(noncritical property from error bounds of solutions under canonical perturbations)\label{theorem2}
Let $(\bar{x},\bar{y})\in S_{\rm KKT}(p_1,p_2)$ with $(p_1,p_2)=(0,0)$ in (\ref{28}).
If there is a number $\varepsilon>0$ and neighborhoods $U$ of $0\in\mathbb{R}^n$,  $W\subseteq S^p$ of 0
such that for any $(p_1,p_2)\in U\times W$ and any $(x_{p_1p_2},y_{p_1p_2})\in S(p_1,p_2)\cap\mathbb{B}_\varepsilon(\bar{x},\bar{y})$ we have the estimate
\begin{equation}\label{101}
\|x_{p_1p_2}-\bar{x}\|=O(\|p_1\|+\|p_2\|)
\end{equation} holds. Then the Lagrange multiplier $\bar{y}\in\Lambda(\bar{x})$
satisfying (\ref{18}) is noncritical for \rm{Problem} (\ref{15}).
\end{theorem}
\noindent{\bf Proof.}
To justify $\bar{y}$ is noncritical, we need to verify by Theorem \ref{theorem1} that the validity of the error bound condition in (\ref{101}) ensures that for any solution pair $(\xi,\eta)\in \mathbb{R}^n\times S^p$ to (\ref{98}) we have $\xi=0$. Pick up any pair ($\xi,\eta)\in \mathbb{R}^n\times S^p$ satisfying (\ref{98}), let $t>0$, and define $(x_t,z_t):=(\bar{x}+t\xi,\bar{z}+tG^\prime(\bar{x})\xi+t\eta)$.
Thus we have for $t$ sufficiently small that
\begin{equation}\label{32}
\begin{array}{lll}
  \Psi_2(x_t,z_t)-\Psi_2(\bar{x},\bar{z})&=&(G(x_t)-\Pi_{S^p_+}(z_t))-(G(\bar{x})-\Pi_{S^p_+}(\bar{z})) \\&=&
tG^\prime(\bar{x})\xi-t\Pi^{\prime}_{S^p_{+}}
 (G(\bar{x})+\bar{y};G^\prime(\bar{x})\xi+\eta)+o(t)\\
 &=&o(t) ~as~t\downarrow 0.
\end{array}
\end{equation}
Since $\Psi_2(\bar{x},\bar{z})=G(\bar{x})-\Pi_{S^p_+}(\bar{z})=0$, we have
\begin{equation}\label{33}
\Psi_2(x_t,z_t)=G(x_t)-\Pi_{S^p_+}(z_t)=-p_{2t} ~with~ p_{2t}=o(t)~ as ~t\downarrow0.
\end{equation}
Moreover,
\begin{equation}\label{34}
\begin{array}{l}
\Psi_1(x_t,z_t)-\Psi_1(\bar{x},\bar{z})\\[4pt]
=
\nabla f(x_t)-\nabla f(\bar{x})+G^{\prime}(x_t)^\ast (z_t-\Pi_{S^p_+}(z_t))-G^{\prime}(\bar{x})^\ast (\bar{z}-\Pi_{S^p_+}(\bar{z}))
\\[4pt]
= t\nabla^2f(\bar{x})\xi+t(G^{\prime\prime}
(\bar{x})\xi)^\ast (G(\bar{x})-\Pi_{S^p_+}(\bar{z}))+G^{\prime}(\bar{x})^\ast ((z_t-\Pi_{S^p_+}(z_t))-(\bar{z}-\Pi_{S^p_+}(\bar{z})))+o(t) \\[4pt]
=t(\nabla^2_{xx}L(\bar{x},\bar{y})\xi+G^\prime(\bar{x})^\ast\eta)
+G^{\prime}(\bar{x})^\ast (tG^{\prime}(\bar{x})^\ast\xi-(\Pi_{S^p_+}(z_t)-\Pi_{S^p_+}(\bar{z})))+o(t)
\\[4pt]
=G^{\prime}(\bar{x})^\ast (tG^{\prime}(\bar{x})^\ast\xi-(\Pi_{S^p_+}(z_t)-\Pi_{S^p_+}(\bar{z})))+o(t)
\\[4pt]
=G^{\prime}(\bar{x})^\ast(tG^{\prime}(\bar{x})^\ast\xi-t\Pi^{\prime}_{S^p_{+}}
 (G(\bar{x})+\bar{y};G^\prime(\bar{x})\xi+\eta))+o(t)=o(t)~ as ~t\downarrow0.
\end{array}
\end{equation}
Since $\Psi_1(\bar{x},\bar{z})=\nabla f(\bar{x})+G^{\prime}(\bar{x})^\ast (\bar{z}-\Pi_{S^p_+}(\bar{z}))=0$, we deduce from the above equality that
\begin{equation}\label{35}
\Psi_1(x_t,z_t)=\nabla f(x_t)+G^{\prime}(x_t)^\ast (z_t-\Pi_{S^p_+}(z_t))=p_{1t}~with~ p_{1t}=o(t)~ as ~t\downarrow0.
\end{equation}
It follow from (\ref{33}) and (\ref{35}) that $(x_t,z_t)$ is a solution to the system $\Psi(x,z)=(p_{1t},-p_{2t})$, where $\Psi(x,z)=(p_{1t},-p_{2t})$ is defined in (\ref{29}). So $(x_t,y_t=z-\Pi_{S^p_{+}}(z))$ is  a solution to (\ref{25}) associated with $(p_{1t},p_{2t})$
and hence we arrive at
\begin{equation}\label{36}
t\|\xi\|=\|x_t\|-\|\bar{x}\|=O(\|p_{1t}\|+\|p_{2t}\|)=
\|o(t)\|
\end{equation}
by (\ref{30}). It yields $\xi=0$ and thus justifies the claim.
\hfill $\Box$

From Theorem \ref{theorem2} we know that if the multiplier for the KKT point is critical then the error bound does not hold at this KKT point. Let us consider the following two examples, and the first one is a convex quadratic semidefinite programming (SDP)  constructed by Bonnans and Shapiro.\begin{example}\label{example2} Consider the following problem [\cite{bonnans2013perturbation}, Example 4.54]:
\begin{equation}\label{102}
\begin{array}{ll}
  \min & x_1+x_1^2+x_2^2 \\
  {\rm s.t.} & {\rm{Diag}(x)}+\varepsilon A\in S^2_{+},
\end{array}
\end{equation}
where $x=(x_1,x_2)\in \mathbb{R}^2$, ${\rm{Diag}}(x)$ is the $2\times2$ diagonal matrix whose ith diagonal element is $x_i$, $i = 1,2$, A is a nondiagonal matrix in $S^2$ , and $\varepsilon$ is a scalar parameter.  When $\varepsilon=0$, the optimization problem (\ref{102}) has the unique optimal solution $x = (0,0)$ with the unique Lagrange multiplier
$\bar{Y}=[\begin{smallmatrix}
     -1 & 0 \\
     0 & 0
\end{smallmatrix}]$. By calculating we get that $\bar{Y}$ is a critical multiplier. It is also easy to see that for any given $\varepsilon\geq0$, problem (\ref{102}) has a unique optimal solution $X(\varepsilon)=(x_1(\varepsilon),x_2(\varepsilon))$ with $x_2(\varepsilon)$ of order $\varepsilon^{2/3}$ as $\varepsilon\rightarrow0$, which illustrates that the error bound like (\ref{101}) does not hold at $(x,\bar{Y})$. So this example is consistent with Theorem \ref{theorem2}.
\end{example}
\begin{example}\label{example3} Consider the following parametric problem:
\begin{equation}\label{103}
\begin{array}{ll}
  \min & x_1^2+x_2^2+x_1x_2-t^{1/2}a\\
  {\rm s.t.} & G(x):=\left[\begin{matrix}
                 x_1^2+x_1x_2 & 0\\
                0 & x_2^2+x_1x_2
               \end{matrix}\right]-tB\in S^2_{+},
\end{array}
\end{equation}
where $x=(x_1,x_2)\in \mathbb{R}^2$, $a:=\frac{5}{3}\sqrt{3} x_1+\frac{4}{3}\sqrt{3}x_2\in\mathbb{R}$, $B=[\begin{smallmatrix}
     2 & 0 \\
     0 & 1
\end{smallmatrix}]\in S^2$ , and $t$ is a scalar parameter.
When $t=0$, the optimization problem (\ref{103}) has a KKT point
$(\bar{x},\bar{Y})=((0,0),O)$ with the multiplier $\bar{Y}=O$ (null matrix).
For any given $t\geq 0$, problem (\ref{103}) has a KKT point $(X(t),Y(t))=((\frac{2}{3}\sqrt{3}t^{1/2},\frac{1}{3}\sqrt{3}t^{1/2}),O)$
which reveals that the error bound like (\ref{101}) holds at $(\bar{x},\bar{Y})$. Calculating based on Corollary \ref{corollary2}, $\bar{Y}$ is a noncritical multiplier. This example is also consistent with Theorem \ref{theorem2}.
\end{example}
In NLPs noncriticality of the multipliers is sufficient for the error bound holding at the KKT point. We inspired by the work of Cui [\cite{Cui2016On}], if we assume some conditions hold, then error bound can be derived from noncritical multipliers.
\begin{theorem}\label{theorem3}
Let $(\bar{x},\bar{y})\in S_{\rm KKT}(p_1,p_2)$ with $(p_1,p_2)=(0,0)$ in (\ref{28}). Define
\begin{equation}\label{43}
\left\{\begin{array}{l}
\mathcal {K}:=(C_{S^n_{+}}(G(\bar{x}),\bar{y}))^{\circ},\\
\Xi:=\{(\xi,\eta)\in \mathbb{R}^n\times S^p :\nabla^2_{xx}L(\bar{x},\bar{y})\xi+G^{\prime}(\bar{x})^{*}\eta=0\},\\
C(\bar{x},\bar{y}):=\{(\xi,\eta)\in \mathbb{R}^n\times S^p :G^{\prime}(\bar{x})\xi\in C_{S^p_+}(G(\bar{x}),\bar{y}), \eta\in C_{S^p_-}(\bar{y},G(\bar{x}))\}.
\end{array}\right.
\end{equation}
 Assume that the following two conditions hold:
 \begin{description}
\setlength{\itemsep}{0pt}
\setlength{\parsep}{0pt}
\setlength{\parskip}{0pt}
\item[(i)] the set $G^{\prime}(\bar{x})^{*}\mathcal {K}$ is closed;\\
\item[(ii)]$\langle\Pi_{\mathcal {K}}(G^{\prime}(\bar{x})\xi),\Pi_{\mathcal {K}}\eta\rangle=0$ for any
$(\xi,\eta)\in\Xi\cap C(\bar{x},\bar{y})$.
\end{description}
If the Lagrange multipliers $\bar{y}\in\Lambda(\bar{x})$
from (\ref{18}) is noncritical for (\ref{15}), then there are numbers $\varepsilon>0$, $\kappa>0$ and neighborhoods $U$ of $0\in\mathbb{R}^n$, $W\subseteq S^p$ of 0
such that for any $(p_1,p_2)\in U\times W$ and any $(x_{p_1p_2},y_{p_1p_2})\in S(p_1,p_2)\cap\mathbb{B}_\varepsilon(\bar{x},\bar{y})$ we have the estimate
\begin{equation}\label{30}
\|x_{p_1p_2}-\bar{x}\|\leq\kappa(\|p_1\|+\|p_2\|).
\end{equation}.\\
In addition, if there exists $\hat{y}\in \Lambda(\bar{x})$ such that $\rank(\bar{x})+\rank(\hat{y})=p$, then we also have
\begin{equation}\label{44}
\dist(y_{p_1p_2}, \Lambda(\bar{x}))=O(\|p_1\|+\|p_2\|).
\end{equation}
\end{theorem}
\noindent{\bf Proof.}
We shall first show that under the given conditions, there exist constants $\varepsilon>0$, $\kappa>0$ and neighborhoods $U$ of $0\in\mathbb{R}^n$ and $W\subseteq S^p$ of 0 such that (\ref{30}) holds.
Suppose on the contrary that (\ref{30}) fails, i.e., for any $k\in\mathbb{N}$ there are $(p_{1k},p_{1k})\in\mathbb{B}_{\frac{1}{k}}(0)\times\mathbb{B}_{\frac{1}{k}}(0_{p\times p})$ and $S_{\rm KKT}(p_{1k},p_{1k})\cap\mathbb{B}_{\frac{1}{k}}(\bar{x},\bar{y})$
satisfying
\begin{equation}\label{54}
\frac{\|x_k-\bar{x}\|}{\|p_{1k}\|+\|p_{2k}\|}\rightarrow\infty
\Leftrightarrow\frac{\|p_{1k}\|+\|p_{2k}\|}{\|x_k-\bar{x}\|}
\rightarrow 0~as~k\rightarrow\infty,
\end{equation}
which yields $p_{1k}=o(\|x_k-\bar{x}\|)$ and $p_{2k}=o(\|x_k-\bar{x}\|)$.
\begin{equation}\label{55}
\begin{array}{lll}
p_{1k}&=&\nabla_xL({x_k},y_k)\\[4pt]
&=&\nabla f(x_k)+G^{\prime}(x_k)^\ast y_k\\[4pt]
&=&\nabla f(\bar{x})+\nabla^2 f(\bar{x})(x_k-\bar{x})
+(G^{\prime}(\bar{x})+G^{\prime\prime}(\bar{x})(x_k-\bar{x}))^\ast
y_k+o(\|x_k-\bar{x}\|)\\[4pt]
&=&\nabla f(\bar{x})+G^{\prime}(\bar{x})^\ast\bar{y}+
G^{\prime}(\bar{x})^\ast(y_k-\bar{y})
+ \nabla^2 f(\bar{x})(x_k-\bar{x})
+G^{\prime\prime}(\bar{x})(x_k-\bar{x}))^\ast
y_k+o(\|x_k-\bar{x}\|)\\[4pt]
&=&\nabla_xL({\bar{x}},\bar{y})+\nabla^2 f(\bar{x})(x_k-\bar{x})
+G^{\prime}(\bar{x})^\ast(y_k-\bar{y})
+(G^{\prime\prime}(\bar{x})(x_k-\bar{x}))^\ast
y_k+o(\|x_k-\bar{x}\|)\\[4pt]
&=&\nabla^2 f(\bar{x})(x_k-\bar{x})+
G^{\prime}(\bar{x})^\ast(y_k-\bar{y})
+(G^{\prime\prime}(\bar{x})(x_k-\bar{x}))^\ast
y_k+o(\|x_k-\bar{x}\|).
\end{array}
\end{equation}
Let us assume without loss of generality that
\begin{equation}\label{60}
\frac{x_k-\bar{x}}{\|x_k-\bar{x}\|}\rightarrow\xi\mbox{ as }k\rightarrow
\infty \mbox{ for some }\xi\neq 0.
\end{equation}
For simplicity, we denote
\begin{equation}\label{56}
t_k:=\|x_k-\bar{x}\|,
\end{equation}
\begin{equation}\label{57}
\Omega:=\{x\in S^p:[P_\beta P_\gamma]^Tx[P_\beta P_\gamma]=0\}
\end{equation}
and for all $k\geq0$,
\begin{equation}\label{58}
\left\{
\begin{array}{ccc}
 G(\bar{x}):=\bar{G}, &G(x_k)+p_{2k}:=G_k, & \tilde{G}_k:=P^TG_kP \\
  \tilde{y}_k:=P^Ty_kP,&H_{1k}:=\Pi_{\Omega}((y_k-\bar{y})/t_k),& H_{2k}:=(y_k-\bar{y})/t_k-H_{1k}\in\mathcal {K}.
\end{array}\right.
\end{equation}
Using Lemma {\ref{Lemma6}} and $\bar{y}\in N(G(\bar{x}),S^p_+)$, ${y_k}\in N(G_k,S^p_+)$ for all $k\geq0$, we deduce that for all $(x_k,y_k)$ sufficiently close to $(\bar{x},\bar{y})$,
\begin{equation}\label{59}
\left\{
\begin{array}{l}
(\tilde{y}_k)_{\alpha\alpha}=O(\|y_k-\bar{y}\|\|G_k-\bar{G}\|), ~ (\tilde{y}_k)_{\alpha\beta}=O(\|y_k-\bar{y}\|\|G_k-\bar{G}\|),  \\
(\tilde{G}_k)_{\beta\gamma}=O(\|y_k-\bar{y}\|\|G_k-\bar{G}\|),
~(\tilde{G}_k)_{\gamma\gamma}=O(\|y_k-\bar{y}\|\|G_k-\bar{G}\|),    \\
(\tilde{y}_k )_{\alpha\gamma}=-\Lambda_\alpha^{-1}(\tilde{G}_k)_{\alpha\gamma}
\Lambda_\gamma+O(\|y_k-\bar{y}\|\|G_k-\bar{G}\|),
\end{array}\right.
\end{equation}
which together with the fact that $(\tilde{G}_k)_{\beta\beta}\in S^{|\beta|}_+$,
$(\tilde{y}_k)_{\beta\beta}\in S^{|\beta|}_-$, yields
\begin{equation}\label{61}
G^{\prime}(\bar{x})\xi \in C_{S^p_+}(G(\bar{x}),\bar{y}),
\end{equation}
\begin{equation}\label{62}
H_1:=\lim_{k\rightarrow\infty}H_{1k}=P\left[\begin{matrix}
                                          0 & 0 & \Lambda_\alpha^{-1}
(P^TG^{\prime}(\bar{x})\xi P) \\
                                          0 & 0 & 0 \\
                                         \Lambda_\alpha^{-1}
(P^TG^{\prime}(\bar{x})\xi P)^T_{\alpha\gamma}\Lambda_\gamma & 0 & 0
                                        \end{matrix}
\right]P^T.\end{equation}
Let us divide the last equality of (\ref{55}) by $\|x_k-\bar{x}\|$ on the both sides. Since
\begin{equation}\label{63}
\nabla^2_{xx}L(\bar{x},\bar{y})\xi=\nabla^2 f(\bar{x})\xi+(G^{\prime\prime}(\bar{x})\xi)^\ast
\bar{y},
\end{equation}
$G^{\prime}(\bar{x})^{\ast}\mathcal {K}$ is assumed to be closed
and that (\ref{61}), (\ref{62}) holds, there exists $H_2\in\mathcal {K}$ such that
\begin{equation}\label{64}
\nabla^2_{xx}L(\bar{x},\bar{y})\xi+G^{\prime}(\bar{x})^{*}(H_1+H_2)=0.
\end{equation}
Let $\eta:=H_1+H_2$. From (\ref{62}) and $H_2\in\mathcal {K}$ we know
\begin{equation}\label{42}
\eta\in C_{S^p_{-}}(\bar{y},G(\bar{x})).
\end{equation}
To verify\begin{equation}\label{66}
(\Sigma_{\alpha\gamma}-E_{\alpha\gamma})\circ(P^T G^\prime(\bar{x})\xi P)_
   {\alpha\gamma}+\Sigma_{\alpha\gamma}\circ(P^T \eta P)_
   {\alpha\gamma}=0
\end{equation} is to verify
\begin{equation}\label{65}
\frac{\lambda_i}{\lambda_i-\lambda_j}\left((P^TG^{\prime}(\bar{x})\xi P)_{ij}-\frac{\lambda_j}{\lambda_i}(P^TG^{\prime}(\bar{x})\xi P)_{ij}\right)=(P^TG^{\prime}(\bar{x})\xi P)_{ij},\forall i\in\alpha, \forall j\in\gamma,
\end{equation}which is obvious.
Then we can obtain from (\ref{61}), (\ref{64}) and (\ref{42}) that $(\xi,\eta)\in\Xi\cap C(\bar{x},\bar{y})$. Furthermore, by using condition (ii) in this theorem we have
\begin{equation}\label{41}
\langle (P^TG^{\prime}(\bar{x})\xi P)_{\beta\beta}, (P^T\eta P)_{\beta\beta}\rangle=\langle\Pi_{\mathcal {K}}(G^{\prime}(\bar{x})\xi ),\Pi_{\mathcal {K}}\eta\rangle=0.
\end{equation}
From (\ref{61}), (\ref{64}), (\ref{42}), (\ref{65}), (\ref{41}) and Theorem
\ref{theorem1} we can get $\xi=0$, which contradicts (\ref{60}). So this justifies
(\ref{30}).\\
Next we shall show that if  there exists $\hat{y}\in \Lambda(\bar{x})$ such that $\rank(\bar{x})+\rank(\hat{y})=p$, then we have (\ref{44}).\\
Denote
\begin{equation}\label{45}
\Delta_1:=\{y\in S^p: \nabla_xL(\bar{x},y)=0\},~
\Delta_2:=\{y\in S^p:y\in N(G(\bar{x});S^p_+)\}.
\end{equation}
Then one has $\Lambda(\bar{x})=\Delta_1\cap\Delta_2$ and $\hat{y}\in\Delta_1\cap \ri(\Delta_2)$. Thus, we obtain from Lemma \ref{Lemma5} that there exists a
constant $\kappa_1>0$ such that
\begin{equation}\label{46}
\dist(y_{p_1p_2},\Lambda(\bar{x}))\leq\kappa_1
(\dist(y_{p_1p_2},\Delta_1)+\dist(y_{p_1p_2},\Delta_2)).
\end{equation}
From Lemma \ref{Lemma3} and the fact that $p_1=\nabla_xL(x_{p_1p_2},y_{p_1p_2})$, we see that
there exist constants $\kappa_2>0$ and $\kappa_2^{\prime}>0$ such that
\begin{equation}\label{47}
\begin{array}{lll}
\dist(y_{p_1p_2},\Delta_1)
&\leq&\kappa_2\|\nabla_xL(\bar{x},y_{p_1p_2})\|\\
&\leq&\kappa_2(\|\nabla_xL(\bar{x},y_{p_1p_2})-
\nabla_xL(x_{p_1p_2},y_{p_1p_2})\|+\|\nabla_xL(x_{p_1p_2},y_{p_1p_2})\|)
\\&\leq&\kappa_2^{\prime}(\|x_{p_1p_2}-\bar{x}\|+\|p_1\|),
\end{array}
\end{equation}
where the last inequality comes from the fact that $f(x)$ and $G(x)$ are both twice continuously differentiable functions. Since $N(\cdot;S^p_-)$ has been proven to be metrically subregular at $\bar{y}$ for $G(\bar{x})$ in Lemma \ref{Lemma8} and $y_{p_1p_2}\in N(G(x_{p_1p_2})+p_2;S^p_+)$, we obtain that there exist constants $\kappa_3>0$ and $\kappa^{\prime}_3\geq0$ such that
\begin{equation}\label{48}
\dist(y_{p_1p_2},\Delta_2)\leq\kappa_3\|G(x_{p_1p_2})+p_2-G(\bar{x})\|
\leq\kappa^{\prime}_3(\|x_{p_1p_2}-\bar{x}\|+\|p_2\|).
\end{equation}
Therefore, combining the inequality (\ref{30}) and inequalities (\ref{46})-(\ref{48}) we show that (\ref{44}) holds. This completes the proof of this theorem.
\hfill $\Box$\\
\begin{Lemma}\label{Lemma4}
Let $\triangle A,~\triangle B,~C\in S^p, ~\bar{A}=\Pi_{S^p_+}(C)$ and $\bar{B}=C-\bar{A}$. Then
$$\triangle A-\Pi^\prime_{S^p_+}(C;\triangle A+\triangle B)=0$$ 
if and only if
\begin{equation}\label{82}
\left\{
\begin{array}{l}
 \triangle A\in C_{S^p_+}(\bar{A}, \bar{B}),\\
  \triangle B+\frac{1}{2}\sigma(\bar{B}, T^2_{S^p_+}(\bar{A};\triangle A))\in[C_{S^p_+}(\bar{A}, \bar{B})]^\circ, \\
 \langle \triangle A, \triangle B\rangle=-\sigma(\bar{B}, T^2_{S^p_+}(\bar{A};\triangle A)).
\end{array}
\right.
\end{equation}
\end{Lemma}
The above Lemma can be obtained from [\cite{ding2017characterization}, Lemma 10] directly.
\begin{theorem}\label{theorem4}
Let $\bar{x}$ be a stationary point of problem (\ref{15}) and let $\bar{y}$ be an associated Lagrange multiplier. If $\bar{x}$ and $\bar{y}$ satisfy the SOSC$\bar{y}$ (\ref{yy}), then the multiplier $\bar{y}$ is noncritical.
\end{theorem}
\noindent{\bf Proof.}
Let $(\xi,\eta)\in\mathbb{R}^n\times S^p$ be arbitrarily chosen such that
\begin{equation}\label{83}
 \left\{\begin{array}{l}
\nabla^2_{xx}L(\bar{x},\bar{y})\xi+G^\prime(\bar{x})^\ast\eta=0, \\
 G^\prime(\bar{x})\xi-\Pi^{\prime}_{S^p_{+}}
 (G(\bar{x})+\bar{y};G^\prime(\bar{x})\xi+\eta)=0.
 \end{array}\right.
\end{equation}
By Lemma \ref{Lemma4}, we know from the first and third lines of (\ref{82})that
\begin{equation}\label{84}
\left\{
\begin{array}{l}
G^\prime(\bar{x})\xi\in C_{S^p_+}(G(\bar{x}),\bar{y}),\\
\langle G^\prime(\bar{x})\xi, \eta\rangle=-\sigma(\bar{y}, T^2_{S^p_+}(G(\bar{x});G^\prime(\bar{x})\xi)).
\end{array}\right.
\end{equation}
Thus, we have $\xi\in C(\bar{x})$. By taking the inner product between $\xi$ and both sides of the first equation of (\ref{83}), respectively, we obtain that
\begin{equation}\label{85}
\begin{array}{lll}
  0&=&\langle \xi,\nabla^2_{xx}L(\bar{x},\bar{y})\xi\rangle+\langle G^{\prime}(\bar{x})\xi,\eta\rangle\\[4pt]
&=&\langle \xi,\nabla^2_{xx}L(\bar{x},\bar{y})\xi\rangle- \sigma(\bar{y}, T^2_{S^p_+}(G(\bar{x});G^\prime(\bar{x})\xi)).
\end{array}
\end{equation}
Hence, it follows from the SOSC$\bar{y}$ (\ref{yy}) for problem (\ref{15}) that $\xi=0$. Therefore the assertion is established.
\hfill $\Box$
\section{A new error bound for $x$-part}
Now, we restrict ourself to the estimation of $x$-part solutions when canonical perturbation is taken.
\begin{theorem}\label{theorem5}
Let $(\bar{x},\bar{y})\in S_{\rm KKT}(p_1,p_2)$ with $(p_1,p_2)=(0,0)$ in (\ref{28}). Then the following properties of problem (\ref{15}) are equivalent:
\begin{description}
\setlength{\itemsep}{0pt}
\setlength{\parsep}{0pt}
\setlength{\parskip}{0pt}
\item[(i)]We can derive $\xi=0$ from
\begin{equation}\label{90}
\left\{\begin{array}{l}
\nabla^2_{xx}L(\bar{x},\bar{y})\xi=0, \\
 G^\prime(\bar{x})\xi-\Pi^{\prime}_{S^p_{+}}
 (G(\bar{x})+\bar{y};G^\prime(\bar{x})\xi)=0.
\end{array}\right.
   \end{equation}
\item[(ii)]There are numbers $\varepsilon>0$, $\ell>0$  and neighborhoods $U$ of $0\in\mathbb{R}^n$, $W\subseteq S^p$ of 0
such that for any $(p_1,p_2)\in U\times W$ and any $(x_{p_1p_2},y_{p_1p_2})\in S(p_1,p_2)\cap\mathbb{B}_\varepsilon(\bar{x},\bar{y})$ we have the estimate
\begin{equation}\label{91}
\|x_{p_1p_2}-\bar{x}\|\leq\ell(\|p_1\|+\|p_2\|+\|y_{p_1p_2}-\bar{y}\|).
\end{equation}
\end{description}
\end{theorem}
\noindent{\bf Proof.}
To verify $(i)\Rightarrow(ii)$, it suffices to check that under the validity of $(i)$, there are numbers  $\varepsilon>0$, $\ell>0$ and neighborhoods $U$ of $0\in\mathbb{R}^n$,  $W\subseteq S^p$ of 0 such that for any $(p_1,p_2)\in U\times W$ and any $(x_{p_1p_2},y_{p_1p_2})\in S(p_1,p_2)\cap\mathbb{B}_\varepsilon(\bar{x},\bar{y})$ we have the estimation (\ref{91}) holds. Assume for the sake of contradiction that for any $k\in\mathbb{N}$ there are $(p_{1k},p_{1k})\in\mathbb{B}_{\frac{1}{k}}(0)\times\mathbb{B}_{\frac{1}{k}}(0_{p\times p})$ and $S_{\rm KKT}(p_{1k},p_{1k})\cap\mathbb{B}_{\frac{1}{k}}(\bar{x},\bar{y})$
satisfying
\begin{equation}\label{92}
\frac{\|x_k-\bar{x}\|}{\|p_{1k}\|+\|p_{2k}\|+\|y_{k}-\bar{y}\|}\rightarrow\infty
\Leftrightarrow\frac{\|p_{1k}\|+\|p_{2k}\|+\|y_{k}-\bar{y}\|}{\|x_k-\bar{x}\|}
\rightarrow 0~as~k\rightarrow\infty,
\end{equation}
which yields $p_{1k}=o(\|x_k-\bar{x}\|)$,  $p_{2k}=o(\|x_k-\bar{x}\|)$
and $\|y_k-\bar{y}\|=o(\|x_k-\bar{x}\|)$.
From (\ref{55}) we have
\begin{equation}\label{93}
p_{1k}=\nabla^2 f(\bar{x})(x_k-\bar{x})+
G^{\prime}(\bar{x})^\ast(y_k-\bar{y})
+(G^{\prime\prime}(\bar{x})(x_k-\bar{x}))^\ast
y_k+o(\|x_k-\bar{x}\|).
\end{equation}
Let us divide the above equality by $\|x_k-\bar{x}\|$ on the both sides and assume without loss of generality that
\begin{equation}\label{94}
\frac{x_k-\bar{x}}{\|x_k-\bar{x}\|}\rightarrow\xi \mbox{ as }k\rightarrow
\infty \mbox{ for some }\xi\neq 0.
\end{equation}
We arrive at
\begin{equation}\label{95}
\nabla^2_{xx}L(\bar{x},\bar{y})\xi=0.
\end{equation}
Since (\ref{59}) and $\|y_k-\bar{y}\|=o(\|x_k-\bar{x}\|)$ we get
\begin{equation}\label{96}
G^\prime(\bar{x})\xi-\Pi^{\prime}_{S^p_{+}}
 (G(\bar{x})+\bar{y};G^\prime(\bar{x})\xi)=0.
 \end{equation}
By assertion (i) of the theorem we know  $\xi=0$
which thus contradicts (\ref{94}). This justifies $(i)\Rightarrow(ii)$.
We next  verify the opposite one $(ii)\Rightarrow(i)$. The proof is similar to that of the Theorem \ref{theorem2}, which just let $\eta=0$.
\hfill $\Box$
\begin{theorem}\label{theorem6}
Let $\bar{x}$ be a stationary point of problem (\ref{15}) and let $\bar{y}$ be an associated Lagrange multiplier. If $\bar{x}$ and $\bar{y}$ satisfy the SOSC$\bar{y}$ (\ref{yy}) then there are numbers $\varepsilon>0$, $\ell>0$  and neighborhoods $U$ of $0\in\mathbb{R}^n$, $W\subseteq S^p$ of 0
such that for any $(p_1,p_2)\in U\times W$ and any $(x_{p_1p_2},y_{p_1p_2})\in S(p_1,p_2)\cap\mathbb{B}_\varepsilon(\bar{x},\bar{y})$ we have the estimate
\begin{equation}\label{97}
\|x_{p_1p_2}-\bar{x}\|\leq\ell(\|p_1\|+\|p_2\|+\|y_{p_1p_2}-\bar{y}\|).
\end{equation}
\end{theorem}
\noindent{\bf Proof.}
From Theorem \ref{theorem4} we have that if $\bar{x}$ and $\bar{y}$ satisfy the SOSC$\bar{y}$ (\ref{yy}), then the multiplier $\bar{y}$ is noncritical.
If the multiplier $\bar{y}$ is noncritical then the assertion (i) in Theorem \ref{theorem5} is satisfied. Using Theorem \ref{theorem5} we have (\ref{97}) which proves this theorem.
\hfill $\Box$
\section{Conclusion}
In this paper we consider a class of generalized SDPs which are not only restricted to linear or convex problems. When the KKT system of the SDP has nonisolated solutions we define the noncritical multipliers and prove under certain conditions, noncriticality of multipliers is equivalent to the local error bound holding at the KKT point, which is useful to construct algorithms that converge locally superlinearly. We introduce a kind of second order sufficient condition SOSC$\bar{y}$ different from the usual SOSC, which is a sufficient condition for noncriticality. Inspired by the structure of noncriticality of multipliers, we show a new structure of KKT points. The equivalence between the new structure and a new error bound for $x$-part has been constructed. This kind of error bound also holds when the KKT point satisfies SOSC$\bar{y}$.\nocite{*}
\bibliographystyle{plain}
\bibliography{refer-1}
\end{document}